\documentclass{article}
\usepackage{graphicx}
\usepackage{color}
\usepackage{amsfonts}
\usepackage{newcent}
\begin{document}
\definecolor{x}{RGB}{10,10,10}
\definecolor{y}{RGB}{250,245,241}
\def \m {\vspace{12pt}

\noindent}
\def \P {{\mathbb P}}
\def \R {{\mathbb R}}
\def \C {{\mathbb C }}
\def \N {{\mathbb N}}
\def \Q {{\mathbb Q}}
\def \Z {{\mathbb Z}}
\def \O {{\cal O}}
\color{x}
\pagecolor{y}

 {\small  
\pagenumbering{roman}
\centerline{Comments about Hilbert's 16'th problem}

\m\m\m\m\m\m\m\m\m
\m \hfill  
$\matrix{
 \hbox{\it John Atwell Moody}
\hfill \cr 
\hbox{\it Warwick Mathematics Institute}
\hfill 
\cr
\hbox{\it  October, 2011}
\hfill
}$
\m
\m\m
\m
\m
\m

}
\vfill\eject

\pagenumbering{arabic}

\m\m{\bf 1. Introduction.}
\m Hilbert's 16'th problem has two familiar  parts
\m{\narrower 
``Was die Curven 6ter Ordnung angeht so habe ich mich -- freilich
auf einem recht umstandlichen Wege..."

}

\m
and
\m{\narrower
   ``Im Anschluss an dieser rien algebraische Problem...''

}

\m
The second part is a hard and unanswered question
of real analysis; it will require a subsequent survey article to this one
to begin to unravel some of what has been done.
\m
The  precursor
was work of Darboux [1] twenty-two
years earlier,  and Poincare [2],
who looked for singular complete curves of degree $e$
invariant under a vector field with poles of degree $d,$
saying that the problem `n'a pas attire l'attention des geometers
autant quelle meritait,' but  `serait
resolu si l'on avait, dans tous les cases, une limite
superior du nombre $e.$' 
\m 
Cerveau and Lins Neto [17] and Carnicer [16] solved Poincare's problem
for nodal curves, and for generic
vector fields. If the curve has nodal singularities, or 
if there are finitely many local invariant 
curves through each singular point of the foliation, $e$ is never more
than\footnote{We are using `degree' to mean the degree of the divisor of poles
of the vector field, which is one less than the number in Carnicer's paper} $d+3.$
Both sets of authors,  I think,  knew that 
$e$ would  not be bounded in terms of $d$ in general.
Lins Neto asked, then [9], 
are
 there are  vector
fields with fixed $d=3,4,...$ and arbitrarily large $e$ 
with no rational first integral. 
\m Lins Neto's question
 was answered
 by  Ollagnier [10], who had
already found  necessary and sufficient
conditions for the special type of quadratic vector field
which are the Lotka-Volterra equations from biology, to have
 no rational first integral, and these include cases when there
are invariant curves of arbitrarily high degree. (An
 independent solution
was is in  [14]).
\m Therefore, even if we restrict the second part of Hilbert's
question to algebraic limit cycles, 
 it is still  not possible to argue that
the absence of any rational first integral for limit cycles
implies them few in number by reasons of degree.
\vfill\eject \noindent  The insightful analysis [11] of Christopher, Llibre and Swirszcz 
  shows that the degree of
 a curve by itself is not a  meaningful invariant after all. The quadratic
transform does not change the dynamical behaviour of
a vector field, while it does change the degree of 
invariant curves. Also the paper includes an interesting
argument over the reals for the existence of an algebraic
curve of degree 12 which is invariant of a particular
quadratic vector field with no rational first integral,
but such that there is no second invariant algebraic curve.

\m 
In this note, we'll 
 restrict to the case of complete invariant
holomorphic curves (not necessarily normal or  irreducible though).
The invariant curve  has a 
 resolution in which it remains invariant and the vector
field does not acquire any new poles, although the foliation does not
become nonsingular. 
\m There is a condition
necessary and sufficient for the existence of a rational first
integral, in terms of one forms on $\P^2$ twisted by $2eH$ for $H$
a line; 
yet it seems difficult to approach this without removing some of the
twisting.
\m
Less extreme log twisting,  once
analytic 
solution germs  $\gamma_i$ are given\footnote{It may be confusing to say they are `given' because then there is at most one solution.  I mean, the 
Zariski closure of the $\gamma_i$ is either $\P^2$ or a lower degree curve,
or else the  $\gamma_i$ can be simultaneously $\delta$ equivariantly deformed.}   at  singular points
$p_1,...,p_s$ of the foliation (here $s\le d^2+3d+3$), shows how to find any  
 complete singular invariant holomorphic curve  of degree 
\footnote{If all $n_i<3$ 
this simplifies to $d+3.$}
 larger
than 
${1\over 2}\sqrt{9+4S
}+d
+{3\over 2}$  passing through only them,  by
rational integration,
where $n_i$ is the maximum by which the local degree of the germ
 exceeds the (positive) order of $\delta$ at points infinitesimally
close to $p_i$
and\footnote{Terms where $n_i\le 0 $ 
should be removed from the sum.} 
\mbox{$S=\sum_{i=1}^s(n_i-1)(n_i-2).$}
\m In  a real polynomial dynamical system of degree d with
a limit cycle, 
if there were ever a  way of learning the maximum 
discrepany $n_1,...,n_s$ over the  set of algebraic
germs, products of at most 
   {\small 
$\left(\matrix{d+2\cr 2}\right)$+1} parts
by Darboux theory,  one could put the largest
$n_1,...,n_s$  in the formula $S.$
As soon as any more than
{\small ${1\over 2}d\sqrt{9+4S}  + {1\over 2}d^2+{1\over 2}S+2$}
limit cycles occur,
 as many as the  remainder  could never be continued to 
complete complex curves, by Harnack's theorem.\footnote{Seven terms cancel.}

\vfill\eject\noindent {\bf 2. Comment on the history of the problem.}
\m
The history is maybe more  important than the question itself. Changes in our interpretation of `integrable' and `rational' 
show the old question in relief. From the first part of Hilbert's question,
about real algebraic curves of degree six,  to second, seems now as if stepping into 
a meaningless abyss, without any reason to have done so. In maths of the Greek empire, people thought that
theorems had political content, and in Hilbert's questions one sees now blind hope for some sort of guidance for
what we should have done, or what will happen next.
\m
A point in the plane, under the influence of a polynomial
vector field, Hilbert must have wondered, should not be just rattling aimlessly,
in the way we now understand the more modern Lorenz
butterfly effect.\footnote{The Poincare-Bendixson theorem
may have been known by then although it was proved  one year later.}
\m
The more modern theory of dynamical systems 
makes a distinction too, though.  It is not any more seeking to find guidance, or analyzing when it may be found or lost, 
but rather accepting and doctrinizing that there is a, possibly natural, system, and within it, points
acting as particles, helplessly obeying it anyway.\footnote{
The Lotka-Volterra model, within which Ollagnier found his counterexample
to the question of Lins Neto, is  this:
there actually is, in the real world, an integer vector made up
of the numbers of organisms of each species in existence. We may
perform a linear regression to see how the logarithmic derivative
 with respect to time may depend on the real value of this vector.
The Lotka-Volterra assumption is that this can be done
without any error, which means that the probability an individual
in one species dying, minus the probability of having a child,
is a utility function  determined  by  linear regression 
on the numbers of individuals of all species. 
}
\m
Now the distinction between separate models is that they each
comprise different hypotheses, differeing from one another as one changes one's mind,  without any continuity; and presenting diverse and discrete choices for human 
intervention into nature, which we may continue to advertise to each other in various ways.

\vfill\eject\noindent
{\bf 3. Chern classes.}
\m If $\delta$ is a meromorphic vector field on the projective
plane with divisor\footnote{Here
$d$ is the order of poles minus the order of zeroes -- yet the `degree'
of the vector field by the usual terminology is  $d+1$ not $d.$} 
$K_1$ of degree  $d,$
the singular subscheme of the
underlying singular foliation is the algebraic cycle

$$ c_2 - (K_1)\cdot (c_1-K_1)$$
where $c_1,c_2$ are  first and second Chern classes 
of the projective plane.  
The three coordinate lines make a triangle of projective
lines, and we can find a one dimensional flow fixing
the corners, or a two dimensional flow preserving the lines
themselves. So we may take as $c_2$ the union of the corners
and as $c_1$ the union of the lines.
The calculation above gives the answer of $c_2$ when $K_1$
is empty and adjusts  it  as $K_1$ enters the picture. 
\m If we take $e$ to  be a natural number, the
necessary and sufficient condition for there to
exist a divisor linearly equivalent to $eH$ with
a component belonging to  a  `level set' of a rational
function, is the condition that there are two
non-linearly-dependent $f,g\in \Gamma(\P^2, \O_{\P^2}(eH))$
such that
$$f\delta(g)=g\delta(f).$$
To make best sense of this we explicitly interpret
global sections of $\O_{\P^2}(eH)$ as rational functions
with poles no worse than $eH.$
\m 
Therefore there is a value of $e$ for which this holds
if and only if there is a rational first integral for
$\delta.$
\m The sheaf of one forms on the projective plane
with poles no worse than $2eH$ has a vector
space of global sections of dimension $4e^2-1,$
and although $fdg$ and $gdf$ may have poles
of order higher than $2eH,$ the difference
$fdg-gdf$ always has poles no worse than $2eH.$
Therefore, we can translate the necessary and sufficient
condition for rational integrability a little bit:
we can find within the projectivication of
$$\Gamma(\P^2,\O_{\P^2}(2eH)\otimes \Omega^1_{\P^2})$$
a copy of the Grassmannian variety of two planes in 
$\Gamma(\P^2,\O_{\P^2}(eH)).$ The two-plane
spanned by $f,g$ is sent to the line spanned by \mbox{$fdg-gdf.$}
\m There is a map of line bundles inducing a map of locally
free coherent sheaves
$$ \Omega_{\P^2} \to \O(K_1)$$
with kernel isomorphic to $\O_{\P^2}(-3H-K_1),$ and
whose cokernel is the actual coordinate ring of the
algebraic cycle  described above (twisted by $K_1$ 
which has a trivial effect).
\m The way that the calculation of the singlar locus 
of the foliation is related to the necessary and
sufficient 
condition for rational integrability 
 is that we  twist
the exact sequence here by $\O_{\P^2}(2eH),$
pass to global sections, and we obtain then  by 
restriction to the Grassmannian subvariety,
a rational map from a Grassmannian variety to
projective space; and the equivalent condition of
existence of a rational first integral is that
this rational map {\it fails} to be a morphism.
\m When there is no rational first integral;
i.e. when the map from the Grassmannian  
is a morphism, the fibers are codimension
zero or dimension zero since the Picard group
of the Grassmannian is cyclic with  (very) ample
generator.\footnote{This observation
was made by D. Maglagan in a seminar, and the subsequent
discussion is based on a further comment by M. Reid; in fact
the comment mentioned Whitney's trick which I now  suspect
could relate
the indeterminacy to higher Chern classes. }
\m There is nothing  interesting to observe from
the standpoint  of dimension only; 
the Grassmannian has dimension $e^2+3e-2,$ 
the projectivized one forms have dimension $4e^2-2,$
and the rational map goes to 
$\P\Gamma\O_{\P^2}(2eH+K_1)$ of dimension {\small$\left(\matrix{2e+d+2\cr
2}\right)$}$-1.$ 
 Although the projectivized one forms twisted in this
manner
do not acquire enough global sections to include the 
Plucker embedding, the projective embedding associated
to the very ample generator of the Picard group, what happens
is that the Grassmannian is still faithfully represented
in the $4e^2-2$  dimensional projective space after a linear
subspace is projected away, perhaps after a projective
normalization. 
\vfill\eject\noindent  {\bf 4. The case of logarithmic twisting}
\m Let us see what can be done if we take global sections
without so much twisting, but only a bit of log twisting.
When we do this, some information will be lost, but not all.
In order to calculate things, we  resolve our curve
let us call it  $C_1$
and so replace $\P^2$ by the surface $S\to \P^2$ obtained
by blowing up singular points of $C_1$ or the total 
transform of $C_1$ which are more than simple crossings.
This does not increase the poles of $\delta,$ in fact since
the chosen points are invariant the flow continues
to preserve the partial resolution.
\m
The total transform of $C_1$ is a simple normal crossing divisor
which we again call $C_1.$ It has coefficients which are positive
numbers;  each prime divisor $P$ in the total
transform of $C_1$ is given the coefficient
which is the valuation at $P$ of any local defining equation of $C_1.$
We will abbreviate this $\nu_P(C_1).$
\m The locally free sheaf of logarithmic one forms 
$\Omega_S(log\  C_1)$ where here $C_1$ abbreviates the pullback 
or total transform of $C_1,$ has a map to $\O(K)$ where $K$ the
divisor of $\delta$ on $S$ is $K_1-Z$ where $Z$ is the divisor of extra
zeroes which $\delta$ acquires on $S.$
\m The first Chern class of $\Omega_S(log\  C_1)$ is the reduced
divisor underlying the total transform of $C_1$ plus the canonical
divisor of the projective plane $-3H$ plus the ramification which is 
the exceptional part of the reduced divisor underlying the total transform.
That is, if we denote the reduced exceptional divisor as $E$
the first chern class is 
$$-3H+E+(C+E)=-3H+E+(eH-I)+E.$$
The first Chern class of the other sheaf is $K=K_1-Z,$ so the
difference, taking $K_1=dH,$ is
$$(e-d-3)H+(2E+Z-I).$$
\vfill\eject\noindent 
The argument of
Bogomolov [6] theorem 12.2 which works in any dimension is also explained in
Miyaoka's paper in a simpler way for surfaces [8]  theorem 2, page
230 attributed to Reid's 
 [7]  Proposition 2, written shortly after his thesis, and attributed
to a discussion\footnote{This refers to the earlier edition of Barth Peters
and Van de Ven.
The current edition  of  Barth Hulek and Van de Ven  includes
now the full text of [7] verbatim.}
 by Barth, Peters and Van de Ven [12] of  Castelnuovo [4] p. 501.  If $w_1,w_2$ are  closed forms in the same line bundle
then  there is a rational function g
with   $w_2=gw_1$ then
 $$0=dw_2=gdw_1+dg\wedge w_1 = dg \wedge w_1 $$
\m
This means that $dg$ is a rational section of the kernel hence
$$0=\delta(g).$$
\m
Although Castelnuovo wrote just five years after Hilbert did,
he likely knew already that Hilbert was asking two unrelated questions.
However, if he had not known that -- if he had wanted to apply his work,
or if he had been motivated
by trying to unite the two parts --  one of the  missing
ingredients would be the degeneration of the Hodge deRham
spectral sequence [5]  3.12 (ii), what is, I think part of Deligne's thesis.
\m
Deligne even observes the specific corollary needed
in ``cas particulier $E^{p0}_1=E^{p0}_\infty$" in Coro 3.2.14, page 39:
logarithmic forms are always closed.
\vfill\eject\noindent 
\m \m To summarize,
\m {\bf 5.  Theorem.}  
For $\delta$  any plane vector field 
of degree $d,$ for  $p_1,...,p_s$ 
singular points of the foliation, and  $\gamma_1,...,\gamma_s$ 
analytic solution germs at $p_1,...,p_s,$
suppose that there is a nontrivial pencil
of plane curves  $C_2$ of degree $q$
whose order at each prime $P$ in the  resolution
of $p_i$ 
is  at least 
$$\nu_P(\gamma_i)-\nu_P(\delta)-2,$$
then any complete solution curve of  degree $q+d+3$ 
passing through only the $\gamma_i$ 
  can be  found by rational integration.
\m{\bf  6. Corollary} For $\delta$  a plane vector field
of degree $d,$ and   $\gamma_1,...,\gamma_s$  analytic
solution germs at $p_1,...,p_s,$ 
let $n_i$ be  the maximum by which   the order of $\gamma_i$
at points  $P$ infinitesimally close to $p_i$ 
exceeds the (positive) order of $\delta.$
Any solution curve  of degree larger \footnote{Terms where $n_i\le 0$
should be removed from the sum.} \footnote{As we remark in another footnote, when all $n_i<3$ this simplifies to $d+3$.} than 
{\small $$
{1\over 2}\sqrt{9+4\sum_{i=1}^s (n_i-1)(n_i-2)
}+d
+{3\over 2}
$$}passing  through only the $\gamma_i$ can be found by rational integration.
\m
Proof. For each number $e$ the complete linear system of degree $e-d-3$
plane curves has dimension {\small  $\left(\matrix{e-d-1\cr 2}\right)-1.$}
The pencil of curves in the theorem includes the pencil of ordinary plane 
curves through the  base cycle  $(n_1-2)p_1+...+(n_s-2)p_s.$
The  dimension of this is  unknown, depending on the special position
of the $p_1,$ but the number is of classical interest. 
Subtracting the sum of the {\small $\left(\matrix{n_1-1\cr 2}\right)$} 
as if the adjunction conditions were linearly independent,
the
familiar quadratic formula gives the value  
of $e$  when the linear
system always begins to move even if the $p_i$ are in general
position.

\m {\bf  7. Remark.} The global  one forms on $\P^2$  with logarithmic
poles on $C_1$ are included
in the global one forms on $S$ with logarithmic poles on $C+E.$
The local monomialization $f=\prod x_i^{e_i}$ 
expresses ${{df}\over f}$ as a linear combination of the ${{dx_i}\over x_i}$
and these patch together.
Therefore any two-plane of global logarithmic forms in the kernel
of $\Omega_{P^2}(log C_1)\to \O(K_1)$
arises geometrically in the manner of the theorem, and all such logarithmic
one forms are closed, even though $\Omega_{\P^2}(log C_1)$ may
not be a locally free sheaf.
\m 

\vfill\eject\noindent
{\bf 8. Remark.} The map $\Omega(log\ C_1)\to \O(K_1)$ exists
even though $C_1$ does not have normal crossings, nor
is $\Omega(log\ C_1)$ locally free. If $F$ is principal
supported on $C_1$ and $f$ the unique nonzero global section
of $\O(F)$ up to scalars then, ${{df}\over f},$ 
maps to a global section of $\O(K_1).$
Darboux observed,  if things are defined over $\Q,$
that once once the number of components of $C_1$ is more than
one (the rank of the Picard)
group plus the dimension of $\Gamma(\O(K_1)),$   there are 
global log forms in the kernel. Also a completely explicit and more
general theorem is in [13]. The kernel is contained in the
kernel of the corresponding map on the resolution and therefore
Darboux' theory can explicitly construct the pencil $C_2$ and 
in turn the indeterminacy on the Grassmannian in cases when 
$C_1$ has many components.

\m {\bf 9. Remark.} It is probably not true that arguments about
global sections of $\Omega(log\ C+E)$ can give the precise 
necessary and sufficient
conditions for a rational first integral in the manner in which
we know that arguments about $\O(2C_1)\otimes \Omega$ can. The technique
of resolutions is just a way of short-cutting the more difficult
Grothendieck group calculation, but the sensible approach may be
to bring in duality and bilinear algebra to try to express the indeterminacy
of the map on the Grassmannian in terms of local contributions where
the vector field meets the curve.

\vfill\eject\noindent \m{\bf 10. Conclusion.}
\m In writing this, I understood a bit about something that
the biologist Jack Cohen once told me. He arrived 
 under the protection of Ian Stewart, and Jack told
me that he had resigned his position 
at another university. His job there, he said, included
assigning a degree classification, a percentage, to each
student. In a meeting, it had been noticed that the distribution
of grades which he assigns does not resemble a Gaussian distribution,
or even a smooth function. He had been asked to assign numbers
a different way.
\m My worry,  along the lines of the 
Lotka-Volterra assumption, had been this: even if we do away
with the linearity assumption,
there is still the issue of {\it Frobenius} integrability
which has to be applied to the construction of the utility
function of each organism, or each person. In the failure
of {\it Frobenius} integrabililty it does not even make
sense to speak of a utility {\it function} at all. 
Jack did not say anything, but he brought me a petri dish,
containing little creatures, hydra, daphnia, and many little
creatures and plants. They were interacting, and almost playing
with each other. 
\m I knew this already, as an American, having come from a country
which was nearly a natural wilderness when I was a child, and I 
had seen even more the overwhelming combination of 
complexity and meaning which is in nature. The indefinite boundary 
between what we thought to be land or sea, 
the creatures, the marshes and rivers, the plants, the tides,
the sunlight,  all changing together in
meaningful ways that could not be described in writing.
\m Jack said, he has been paid, and taken money, to help
create duck ponds, where a stream passes through. He said,
there is usually a sort of machine, so that the water has to
go underground to leave the lake, and owners hate this, they
remove it. But then the trough beneath begins to fill with
silt, and they call on Jack as a consultant. He said, he
takes the payment and advises them to reinstate the machinery.
He said, there is such a duck pond here; it was the pride of our
university. Since then a small area
has also returned to nature.

\vfill \eject
\noindent \pagestyle{empty}\normalsize
[1] G. Darboux, Memoire sur les equations di.erentielles algebriques du premier
ordre et du premier degre (Melanges) , Bull. Sci. math. 2`eme serie 2 (1878),
60.96; 123.144; 151.200

\vspace{6pt}
\noindent 
[2]  H. Poincare, Sur l.integration algebrique des equations differentielles,
Comptes Rendus des Seances, Academie des Sciences, Monday, 15 April, 1891.

\vspace{6pt}
\noindent
[3] D. Hilbert, Paris ICM, 1900 (from wikipedia)

\vspace{6pt}
\noindent
[4] G. Castelnuovo, Sulle Superficie aventi il genere aritmetico negative, Rend Circ Mat Palermo 20 (1905) Memorie Scelte \#XXVII

\vspace{6pt}
\noindent
[5] Deligne, Theorie de Hodge II, Publications mathematiques de l'IHES, 40 (1971)

\vspace{6pt}
\noindent
[6] Bogomolov, Holomorphic tensors and vector bundles, Math USSR Izvestija 13 (1979)

\vspace{6pt}
\noindent 
[7] Reid, Bogomolov's theorem $c_1^2<4c_2$, Int's Symp on  Algebraic Geometry, Kyoto (1977)

\vspace{6pt}
\noindent
[8] Miyaoka, On the Chern numbers of surfaces of general type, Inventiones 42 (1977)

\vspace{6pt}
\noindent 
[9]  A. Lins Neto, Some examples for Poincare and Painleve problems, Ann. Scient. Ec. Norm. Sup. 35 (2002), 231.266

\vspace{6pt}
\noindent 
[10] J.M. Ollagnier, About a conjecture on quadratic vector fields,
Journal of Pure and Applied Algebra 165 (2001)

\vspace{6pt}
\noindent
[11] C. Chrisotpher, J. Llibre, G. Swirszcz, Invariant curves of large degree for quadratic system,
Math Analysis and Applications 202 (2005), 

\vspace{6pt}
\noindent
[12]  Barth, Peters, Van de Ven, compact complex surfaces (1984)

\vspace{6pt}
\noindent
[13] J. Llibre, X. Zhang, Rational first integrals in the Darboux theory of integrability in Cn, Bulletin des Sciences Mathematiques
Volume 134, Issue 2, March 2010, Pages 189-195

\vspace{6pt}
\noindent
[14] C. Christopher and J. Llibre, 
A family of quadratic polynomial differential systems 
with invariant algebraic curves of arbitrarily high degree without rational
first integrals, Proc. Amer. Math. Soc., 130(2002), 2025-2030.

\vspace{6pt}
\noindent 
[15] Camacho, Cesar; Sad, Paulo
Invariant varieties through singularities of holomorphic vector fields. 
Ann. of Math. (2) 115 (1982), no. 3, 579.595. 

\vspace{6pt}
\noindent
[16] M. Carnicer,  The Poincare problem
in the nondicritical case, Annals of Mathematics 1994.

\vspace{6pt}
\noindent
[17] D. Cerveau, A. Lins Neto, Holomorphic foliations in CP2
having an invariant algebraic curve, Ann. Inst. Fourier, Grenoble 41 (1991)
\end{document}